\documentclass[twoside, 11pt]{article}
\usepackage{amsmath} 
\usepackage[e]{esvect}
\usepackage{graphicx}
\usepackage{epstopdf}
\usepackage[adobe-utopia]{mathdesign}
\usepackage[T1]{fontenc}
\usepackage[english]{babel}
\usepackage{color}

\usepackage{hyperref}
\usepackage[a4paper, total={5in, 9in}]{geometry}
\title{\textsc{A note on a flip-connected class of generalized domino tilings of the box $[0,2]^n$}}
\date{}
\author{{Andrzej P. Kisielewicz}\\
{\small \sf  A.Kisielewicz@wmie.uz.zgora.pl}\\
{\small\it Wydzia{\l} Matematyki, Informatyki i Ekonometrii, Uniwersytet Zielonog\'orski,}\\
{\small\it ul. Podg\'orna 50, 65-246 Zielona G\'ora, Poland}}

\newtheorem{lemat}{\sc Lemma}
\newtheorem{tw}{\sc Theorem}

\newtheorem{wn}{\sc Corollary}
\newtheorem{df}{\sc Definition}

\newtheorem{uw}{\sc Remark}
\newtheorem{uwi}{\sc Remarks}

\newtheorem{nap}{\sc Example}

\newtheorem{nps}[nap]{\sc Examples}

\def\ka #1{\mathscr{#1}}
\def\kal #1 #2{\mathscr{#1}^{#2}}
\def\proof{\noindent \textit{Proof.\,\,\,}}

\def\enn{\mathbb{N}}

\def\er{\mathbb{R}}

\def\Aut #1 #2{\operatorname{Aut}^{#1} (#2)}

\def\prop #1{\operatorname{prop}\, #1}

\def\bred #1 {\colorbox{red}{ #1}}
\def\red #1 {{\color{red} #1 }} 

\graphicspath{{figpdfkolor/}}
\DeclareGraphicsExtensions{.pdf,.png} 
\begin{document}
\maketitle
%\begin{frontmatter}
%title{\textsc{There are at most $2^{d+1}-2$ neighbourly simplices in dimension $d$}}
%\date{}
%\author{{Andrzej P. Kisielewicz \,\, Krzysztof Przes{\l}awski}\\
%{\small\it Wydzia{\l} Matematyki, Informatyki i Ekonometrii, Uniwersytet Zielonog\'orski,}\\
%{\small\it ul. Podg\'orna 50, 65-246 Zielona G\'ora, Poland}\\
%{\small \sf  A.Kisielewicz@wmie.uz.zgora.pl K.Przeslawski@wmie.uz.zgora.pl}}

%\end{frontmatter}

\begin{abstract}
Let $n,d\in \enn$ and $n>d$. An $(n-d)$-domino is a box $I_1\times \cdots \times I_n$ such that $I_j\in \{[0,1],[1,2]\}$ for all $j\in N\subset [n]$ with $|N|=d$ and $I_i=[0,2]$ for every $i\in [n]\setminus N$. If $A$ and $B$ are two $(n-d)$-dominoes such that $A\cup B$ is an $(n-(d-1))$-domino, then $A,B$ is called a twin pair. If $C,D$ are two $(n-d)$-dominoes which form  a twin pair such that $A\cup B=C\cup D$ and $\{C,D\}\neq \{A,B\}$, then the pair $C,D$ is called a flip of $A,B$. A family $\ka D$ of $(n-d)$-dominoes is a tiling of the box $[0,2]^n$ if interiors of every two members of $\ka D$ are disjoint and $\bigcup_{B\in \ka D}B=[0,2]^n$. An $(n-d)$-domino tiling $\ka D'$ is obtained from an $(n-d)$-domino tiling $\ka D$ by a flip, if there is a twin pair $A,B\in \ka D$ such that $\ka D'=(\ka D\setminus \{A,B\})\cup \{C,D\}$, where $C,D$ is a flip of $A,B$. A family of $(n-d)$-domino tilings of the box $[0,2]^n$ is flip-connected, if for every two members $\ka D,\ka E$ of this family the tiling $\ka E$ can be obtained from $\ka D$ by a sequence of flips. In the paper some flip-connected class of $(n-d)$-domino tilings of the box $[0,2]^n$ is described.
%Based on local transformations of some class of partitions of the box $[0,2]^d$ into boxes, where $d\leq 6$, some flip-connected class of generalized domino tilings of $[0,2]^n$, where $n>d$, is described. 
\end{abstract}

\section{Introduction}
A non-empty set of the form $K_1\times\cdots \times K_d \subset [0,2]^d$, where $K_i\subset [0,2]$ for $i\in [d]=\{1,...,d\}$, is called a {\it box}. A {\it domino} is a box $I_1\times \cdots \times I_n\subset [0,2]^n$ such that $I_i=[0,2]$ for some $i\in [n]$ and $I_j\in \{[0,1],[1,2]\}$ for all $j\in [n]\setminus \{i\}$. A natural generalization of a domino is an $(n-d)$-{\it domino}, where $n>d$, $n,d\in \enn$, which is a box $I_1\times \cdots \times I_n$ such that $I_j\in \{[0,1],[1,2]\}$ for all $j\in N\subset [n]$ with $|N|=d$ and $I_i=[0,2]$ for every $i\in [n]\setminus N$. Thus, a $1$-domino is a domino. Let $A$ and $B$ be two $(n-d)$-dominoes such that $A\cup B$ is an $(n-(d-1))$-domino. In this case there is $i\in [n]$ such that $A_j=B_j$ for all $j\in [n]\setminus \{i\}$ and $A_i,B_i\in \{[0,1],[1,2]\}$, $A_i\neq B_i$. Such pair $A,B$  will be called a {\it twin pair} (Figure 1). This $(n-(d-1))$-domino $A\cup B$ can be divided into two new $(n-d)$-dominoes $C$ and $D$ (Figure 1), that is, $A\cup B=C\cup D$,  $C_j=D_j$ for all $j\in [n]\setminus \{k\}$ and $C_k,D_k\in \{[0,1],[1,2]\}$, $C_k\neq D_k$, where $k\neq i$. The twin pair $C,D$  is called a {\it flip} of $(n-d)$-dominoes $A$ and $B$. Let us observe that there are $(n-d)$ flips of $A,B$ (Figure 1). 

\vspace{-1mm}
{\center
\includegraphics[width=7cm]{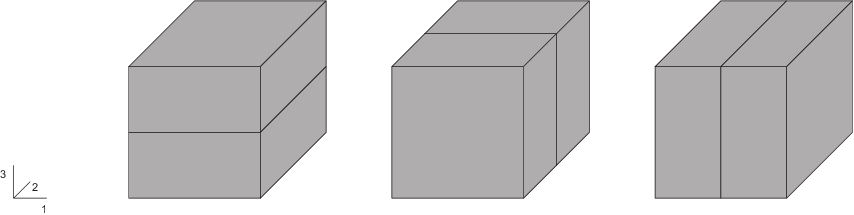}\\
}

\smallskip
\noindent{\footnotesize Figure 1: Two $2$-dominoes $A=[0,2]^2\times [0,1]$ and  $B=[0,2]^2\times [1,2]$ which form a twin pair (on the left) and two flips of $A,B$: $C=[0,2]\times [0,1]\times [0,2], D=[0,2]\times [1,2]\times [0,2]$ (in the middle) and $E=[0,1]\times [0,2]^2, F=[1,2]\times [0,2]^2$ (on the right).}

\medskip
A family $\ka D$ of $(n-d)$-dominoes is a {\it tiling} of the box $[0,2]^n$ if interiors of every two members of $\ka D$ are disjoint and $\bigcup_{B\in \ka D}B=[0,2]^n$. Clearly, $|\ka D|=2^d$.  Let us suppose that an $(n-d)$-domino tiling $\ka D$ contains a twin pair $A,B$. Flipping this pair we obtain a new $(n-d)$-domino tiling $\ka D'$ of $[0,2]^n$, where $\ka D'=(\ka D\setminus \{A,B\})\cup \{C,D\}$ and $C,D$ is a flip of $A,B$. %Now, having a second $(n-d)$-domino  tiling  $\ka D'$ we may try to pass by flipping from $\ka D$ to $\ka D'$.  
The aim of this note is to present some class of $(n-d)$-domino tilings of the box $[0,2]^n$ which is {\it flip-connected}, that is, for every two tilings $\ka D, \ka E$ of this class there is a sequence $\ka D=\ka T_1,...,\ka T_m=\ka E$ of  $(n-d)$-domino tilings of the class and there are two sequences $A_1,B_1,....,A_{m-1},B_{m-1}$ and $C_1,D_1,...,C_{m-1},D_{m-1}$ of twin pairs, where $A_i,B_i\in \ka T_i$ for $i\in [m-1]$, such that the pair $C_{i},D_{i}$ is a flip of the pair $A_i,B_i$ for $i\in [m-1]$, and moreover $\ka T_{i+1}=(\ka T_i\setminus \{A_i,B_i\})\cup \{C_{i},D_{i}\}$ for $i\in [m-1]$. 

%of flips which transforms one into the another (Corollary \ref{alm}). 

\smallskip
There are various tilings which can be modify by a sequence of local transformations. The most known are domino and lozenge tilings (\cite{R,W}). Local transformations of  domino tilings in dimensions $d\geq 3$ are also examined (\cite{F,KS,MS,Sa,ST}). Generalized domino tilings are considered in \cite{Pa}.

\smallskip
Our result for $(n-d)$-dominoes is based on special codes which encode some partitions of a $d$-dimensional box $X$ into smaller boxes  (\cite{LS2,KP}).
Let $S$ be a set of arbitrary objects, and let $s\mapsto s'$ be a permutation  of $S$ such that $s''=(s')'=s$ and $s'\neq  s$. Let $S^d=\{v_1...v_d\colon v_i\in S\}$.  Elements of $S$ are called {\it letters}, while members  of $S^d$ will be called {\it words}. We add to $S$ an extra letter $*$ and the set $S\cup \{*\}$ will be denoted by $*S$. We assume that $*'=*$. (Clearly, elements of $(*S)^d$ will be also called words.) Two words $v,w\in (*S)^d$ are {\it dichotomous} if there is $i\in [d]$ such that $v_i,w_i\in S$ and $v_i=w'_i$. If additionally $v_j=w_j$ for $j\in [d]\setminus \{i\}$, then we say that $v$ and $u$ form a {\it twin pair} ({\it in the $i$-th direction}). A set $V\subset (*S)^d$ is a {\it code} if every two distinct words in $V$ are dichotomous (compare Table 1). It is easy to see that $|V|\leq 2^d$. If $V\subset S^d$ and $|V|=2^d$, then $V$ is called a {\it cube tiling} code. If $V$ is a cube tiling code such that $v_i\in \{w_i,w_i'\}$ for every $v,w\in V$ and every $i\in [d]$, then $V$ is called {\it simple}.

A $d$-{\it box} is a set of the form $X=X_1\times \cdots \times X_d$, where $|X_i|\geq 2$ for $i\in [d]$. A {\it box} in a $d$-box $X$ is any non-empty set  of the form $K=K_1\times\cdots \times K_d \subset X$, where $K_i\subset X_i$ for $i\in [d]$. A box $K$ is \textit{ proper} if $K_i\neq X_i$ for each $i\in [d]$. Two boxes $K, G\subset X$ are called \textit{dichotomous} if there is $i\in [d]$ such that $K_i=X_i\setminus G_i$. If $\ka F$ is a family of pairwise dichotomous proper boxes in $X$ containing $2^d$ boxes, then it is called a {\it minimal partition} of $X$ (Figure 2). Each minimal  partition of $X$ is a partition of the $d$-box $X$, that is, $\bigcup_{K\in \ka F}K=X$. 
 
Cube tiling codes encode minimal partitions. To see this, let $\ka P(A)$ be the power set of a set $A$, and let for each $i\in[d]$ a mapping $f_i\colon S\to \ka P(X_i)\setminus \{\emptyset,X_i\}$ be such that $f_i(s')=X_i\setminus f_i(s)$. %Additionally, $f_i(*)=X_i$ for $i\in [d]$. 
We define the mapping $f\colon S^d\to \ka P(X)$ by $f(s_1\ldots s_d)=f_1(s_1)\times\cdots\times f_d(s_d).$ The set $f(V)=\{f(v)\colon v\in V\}$ is said to be a \textit{realization} of the code $V$ (Figure 2 and Table 1). If $V\subset S^d$ is a cube tiling code, then the set of boxes $f(V)$ is a minimal partition of $X$. Cube tiling codes encode also other interesting combinatorial structures such as  $2$-periodic cube tilings of $\er^d$ or perfect codes in the maximum metric (\cite{Co,Kii}). 

It is easy to check that any minimal partition defines a cube tiling code (see \cite[Subsection 2.6]{Ki2}).

Two dichotomous boxes $K, G\subset X$ are called a \textit{twin pair} if there is $i\in [d]$ such that $K_i=X_i\setminus G_i$ and  $K_j=G_j$ for all $j\in [d]\setminus \{i\}$.  Let us note that for any twin pair $K,G$, where  $K_i=X_i\setminus G_i$, the set $H=K\cup G$ is a box with $H_i=X_i$. The box $H$ is called a {\it gluing} of $K,G$. Taking a non-empty set $M_i\subset X_i$, where $M_i\not\in \{X_i, K_i, G_i\}$  we can divide the box  $H$ into two new boxes $M$ and $N$, that is $H=M\cup N$, where $M_j=N_j=K_j$ for $j\in [d]\setminus \{i\}$ and $N_i=X_i\setminus M_i$. Let us call this operation a {\it cutting} of $H$. Thus, if $\ka F$ is a minimal partition of a $d$-box $X$ and a twin pair $K,G$ belongs to $\ka F$, then by gluing and cutting we can pass from $\ka F$ to a new minimal partition $\ka G=(\ka F\setminus \{K,G\})\cup \{M,N\}$ (\cite{Kgl}).  In Figure 2 we see a passing by gluing and cutting from one minimal partition to another.

\vspace{-2mm}
{\center
\includegraphics[width=10cm]{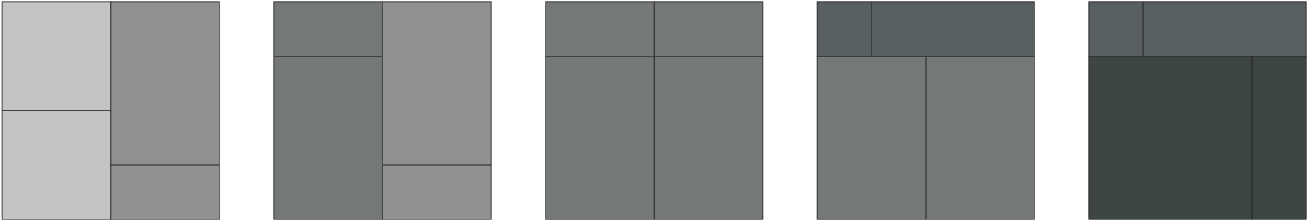}\\
}

\medskip
\noindent{\footnotesize Figure 2: A passing by gluing and cutting from one minimal partition of the box $[0,2]^2$ to another. At each step of the process precisely one twin pair is modified. The partitions on the most left and the most right are a realization of the codes $V=\{aa,aa',a'b,a'b'\}$ and  $W=\{cc, c'c, bc', b'c'\}$, respectively. In the realizations of the codes we took $f_1(a)=f_2(a)=[0,1], f_1(b)=f_2(b)=[0,\frac{1}{2}]$ and $f_1(c)=f_2(c)=[0,\frac{3}{2}].$}

\smallskip
%Clearly, the procedure of gluing and cutting can be expressed in the language of codes ([]). 
For simplicity we shall work with codes $V\subset S^d$ rather than boxes. Moreover, if $v,w$ is a twin pair with $v_i=w_i'$, then a twin pair $u,q$ such that $u_i=q_i'$,  $\{v_i,w_i\}\neq \{u_i,q_i\}$ and $v_j=w_j=u_j=q_j$ for $j\in [d]\setminus \{i\}$ will be called a {\it shift} of the pair $v,w$ (\cite{Ost}). If $u,q$ is a shift of $v,w$, then we will write $\{v, w\}\rightarrow \{u, q\}$. Thus, in some cases we may pass from a cube tiling code $V$ to a cube tiling code $W$ by a sequence of shifts. For example, we may pass from $V$ to $W$ (Figure 2) by shifts:

\vspace{-4mm}
$$
\{aa, aa'\}\rightarrow \{ac, ac'\},\;\; \{a'b, a'b'\}\rightarrow \{a'c, a'c'\},\;\; \{ac', a'c'\}\rightarrow \{bc', b'c'\}, \;\; \{ac, a'c\}\rightarrow \{cc, c'c\}.
$$ 

If $W$ is obtained from $V$ by shifts, then there are two sequences of twin pairs $(u^n,q^n)_{n=1}^k$, $(v^n,w^n)_{n=1}^k$ and there is a sequence of cube tiling codes $(V^j)_{j=1}^{k+1}$ such that:

\smallskip
$\bullet$  For every $j\in [k]$ the twin pair $v^j,w^j$ is a shift of the twin pair $u^j,q^j$.  

\smallskip
$\bullet$ $V^1=V$, $V^{j+1}=(V^{j}\setminus \{u^j,q^j\})\cup \{v^j,w^j\}$, where $\{u^j,q^j\}\subset V^{j}$ for every $j\in [k]$ and $V^{k+1}=W$.
%and a sequence of cube tiling codes $(V^j)_{j=1}^{k+1}$ such that 

\section{Flips of regular $(n-d)$-domino tilings}

% and we can do that not only for cube tiling codes but also for smaller codes. More pre and cisely, two codes $V,W\subset S^d$ are {\it equivalent} if $\bigcup f(V)=\bigcup f(W)$ for every above defined mapping $f$.
\vspace{-2mm} 
In  \cite[Theorem 1.5]{Kgl} (see also \cite{DIP,Ost}) we proved the following

\vspace{-2mm}
\begin{tw}
\label{glue}
Let $d\leq 6$. For every cube tiling codes $V,W\subset S^d$ one may pass by shifts from $V$ to $W$.
\end{tw}

We shall adapt this theorem to the case of $(n-d)$-domino tilings of $[0,2]^n$.

\vspace{-2mm}
\begin{uw}{\rm In \cite{Kgl} we use a slightly different terminology and if it is possible to pass by shifts from $V$ to $W$, then we say that $V$ and $W$ are {\it strongly equivalent.}
}

\end{uw}

Let $*A=\{0,1,*\}$, and  let $v\in  (*A)^n$. A geometric interpretation of the word $v=v_1...v_d$ is the box $\breve{v}=\breve{v}_1\times \cdots \times \breve{v}_n\subset [0,2]^n$, where $\breve{v}_i=[0,1]$ if $v_i=0$, $\breve{v}_i=[1,2]$ if $v_i=1$ and $\breve{v}_i=[0,2]$  if $v_i=*$.  Let $\prop (v)=\{i\in [n]\colon v_i\neq *\}$ for $v\in (*A)^n$.  If $|\prop (v)|=d$, then the box $\breve{v}$ is an $(n-d)$-domino. To make notations shorter, we shall be  working mainly with $(n-d)$-dominoes in the form of words rather than boxes. Thus, a family of $(n-d)$-dominoes in which every pair of different members have disjoint interiors   will be represented by a (dichotomous) code $D\subset (*A)^n$, where $d=|\prop (v)|$ for every $v\in D$ and $0'=1$ (compare Table 1). 

%\begin{uw} {\rm
%In \cite{Als} (see also \cite{Pr}) it was shown that if $D\subset (*A)^n$ is an $(n-d)$-domino tiling of $[0,2]^n$ and $\prop(D)=\bigcup_{v\in D}\prop(v)$, then $|\prop(D)|\leq 2^d-1$ and for every $m\in \{d,..., 2^{d}-1\}$ there is an  $(n-d)$-domino tiling $D$ such that $|\prop(D)|=m$.
%}
%\end{uw}

Let  $D$ be an $(n-d)$-domino tiling of the box $[0,2]^n$, and let $M(D)=[d_{ij}]_{i\in [2^d], j\in [n]}$ be a matrix in which  rows are words from $D$. We call the tiling $D$ {\it regular} (see Example 2 and Table 1) if there is a partition of the set $[n]$

\vspace{-3mm}
\begin{equation}
\label{roz}
[n]=N_1\dot{\cup} \cdots \dot{\cup} N_d \dot{\cup} N^*
\end{equation}

with the following properties:

\medskip
($\alpha$) For every $i\in [d]$ the set $N_i=\{j_1<...<j_{n_i}\}$ is non-empty, and the set $K_i$ of the rows of the sub-matrix $[d_{ij}]_{i\in [2^d], j\in N_i}$ of the matrix $M(D)$ %consisting of all $j$-th columns of $M(D)$, where  $j\in N_i=\{j_1<...<j_{n_i}\}$,  
is of the form: 

\vspace{-4mm}
\begin{equation}
\label{ki}
K_i=\{\delta *\ldots *,\;\; *\delta *\ldots *, ..., *...*\delta \colon \delta \in \{0,1\}\}.
\end{equation}

\medskip
($\beta$) The set $N^*$ consists of positions of $*$-columns in $M(D)$ and if $N^*\neq\emptyset$, then the set $K^*$ of the rows of the sub-matrix $[d_{ij}]_{i\in [2^d], j\in N^*}$ of the matrix $M(D)$ contains one word $*...*$ of length $|N^*|$.

%Moreover, if $N^*\neq\emptyset$, then the set $K^*$ of the rows of the sub-matrix $M(D)_{2^d\times |N^*|}$ contains one word $*...*$ (see Table 1).  

\medskip
The enumeration of the sets $N_1,...,N_d$ is such that $1=i_1<...<i_d$, where $i_1,...,i_d$ are the first elements in $N_1,...,N_d$.  Let us note that $|\prop(u)\cap N_i|=1$ for every  $u\in D$ and $i\in [d]$, $|K_i|\geq 2$ for $i\in [d]$ (if $|K_i|=2$, then $K_i=\{0,1\}$)  and $|K^*|\geq 0$. 

\vspace{-2mm}
\begin{nap}{\rm Since it can happen that $K_i=\{0,1\}$ for some $i\in [d]$ and $K^*\neq\emptyset$, it follows that  for every $n,d\in \enn$, $n>d$, there is a regular $(n-d)$-domino tiling of $[0,2]^n$, as  $D=\{\varepsilon_1...\varepsilon_d*^{n-d}\colon (\varepsilon_1,...,\varepsilon_d)\in \{0,1\}^d\}$ is regular. In this case we have $N_i=\{i\}$, $K_i=\{0,1\}$ for $i\in [d]$ and $N^*=\{d+1,...,n\}$, $K^*=\{*^{n-d}\}$. Let us note that having a regular $((n-1)-(d-1))$-domino tiling $D$ of the box $[0,2]^{n-1}$ one may easily construct a regular $(n-d)$-domino tiling $E$ of the box $[0,2]^{n}$. To do that, let $[n-1]=M_1\dot{\cup} \cdots \dot{\cup} M_{d-1} \dot{\cup} M^*$ be the partition (\ref{roz}) for the tiling $D$, and let $K_1,...,K_{d-1}$ be the sets of the form (\ref{ki}) for $D$. Then taking the partition $[n]=N_1\dot{\cup} \cdots \dot{\cup} N_d \dot{\cup} N^*$, where $N_1=\{1\}$, $N_{i+1}=M_i+1$ for $i\in [d-1]$ and $N^*=M^*+1$ ($A+1=\{a+1\colon a \in A\}$) we obtain the regular  $(n-d)$-domino tiling $E$ of the box $[0,2]^n$, where sets (\ref{ki})  are of the form $G_1=\{0,1\}, G_{i+1}=K_i$ for $i\in [d-1]$ and $G^*=K^*$.

The above constructed tilings have the property that $K_i=\{0,1\}$ for some $i\in [d]$. It is easy to check that for $n=2,3$ all $(n-d)$-domino tilings of $[0,2]^n$ have that property. For example, in the regular domino tiling $D=\{0*0,1*0,*01,*11\}$ of $[0,1]^3$  we have $K_1=\{0*,1*,*0,*1\},K_2=\{0,1\}$ ($N_1=\{1,2\},N_2=\{3\}$). The first dimension $n\in \enn$ in which there is an $(n-d)$-domino tiling $D$ of $[0,2]^n$ such that the matrix $M(D)$ does not contain a column which  contains only $0,1$ is $n=4$ (compare Example 2). For example, for $n=4$ and $d=3$ the set $D=\{*01\varepsilon,1*0\varepsilon,01*\varepsilon\colon \varepsilon \in \{0,1\}\}\cup \{000*,111*\}$ is a domino tiling of $[0,2]^4$ and there is no column in the matrix $M(D)$ which contains only 0 and 1. Let us note that the tiling $D$ is not regular. 
}
\end{nap}

If the matrix $M(D)$ of an $(n-d)$-domino tiling $D$ of $[0,2]^n$ contains a $*$-column, then one can obtain from the tiling $D$ a new $((n-1)-d)$-domino tiling $E$ of $[0,2]^{n-1}$, where $E=\{u_1...u_{j-1}u_{j+1}...u_n\colon u\in D\}$ and $j\in N^*$ (that is, $u_j=*$ for every $u\in D$). Every $(n-d)$-domino tiling $D$ of $[0,2]^n$ such that $M(D)$ contains a $*$-column will be called {\it $*$-reducible}, and otherwise $D$ will be called {\it $*$-irreducible}.

As we have just seen, since we admit that $K^*\neq\emptyset$, we can construct  a regular $(n-d)$-domino tiling of $[0,2]^n$ for every $n,d\in \enn$, $n>d$.  
This is no longer true for regular $(n-d)$-domino tilings of $[0,2]^n$ which are $*$-irreducible, that is, $K^*=\emptyset$. In the next section we give a characterization of numbers $n,d\in \enn$, $n>d$, for which regular $*$-irreducible $(n-d)$-domino tilings of the box $[0,2]^n$ exist.

\subsection{Regular $(n-d)$-domino tilings and cube tiling codes}

%In this section we describe the relationship between regular $(n-d)$-domino tilings of the box $[0,2]^n$ and cube tiling codes $V\subset S^d$.
%In this section, based on the regular structure of the sets $K_1,...,K_d$, we describe the relationship between regular $(n-d)$-domino tilings of the box $[0,2]^n$ and cube tiling codes $V\subset S^d$.
The structure of sets $K_1,...,K_d$ makes regular $(n-d)$-domino tilings of the box $[0,2]^n$ very close to cube tiling codes $V\subset S^d$. In this section we describe the relationship between these two types of tilings.

Let $S=\{a_1,a'_1,...,a_n,a_n'\}$, and let $D$ be a regular $(n-d)$-domino tiling of the box $[0,2]^n$. For every $u=u_1...u_n\in D$ we define  a word $v(u)=v_1...v_d\in S^d$: For each $i\in [d]$ let

\vspace{-2mm}
\begin{equation}
\label{rel}
v_i=
\begin{cases}
a_j & \text{if $j\in N_i$ and $u_j=0$},\\
a_j' & \text{if $j\in N_i$ and $u_j=1$},\\
\end{cases}
\end{equation}
where $N_i$, $i\in [d]$, is such in (\ref{roz}).  It is easy to show that the set $V_D=\{v(u)\colon u\in D\}$ is a cube tiling code.To do that, let $u,q\in D$ be two distinct words. Then there are $i\in [d]$ and $j\in N_i$ such that $u_j=0$ and $q_j=1$, and consequently $v_i(u)=v_i(q)'$, that is, the words $v(u)$ and $v(q)$ are dichotomous, and since $|D|=2^d$, we have $|V_D|=2^d$. Thus, $V_D$ is a cube tiling code. We shall say that $V_D$ is the {\it cube tiling code of } the tiling $D$.

\vspace{-3mm}
\begin{nap} {\rm
For the code $D$ with words given in the matrix $M(D)$ on the left in Table 1 we have $N_1=\{1,2\}, N_2=\{3,4\}$, $N_3=\{5,6\}$ and $N^*=\{7\}$. Moreover, $K_i=\{0*,1*,*0,*1\}$ for $i\in \{1,2,3\}$ and $K^*=\{*\}$. If $S=\{a_1,a_1',...,a_7,a_7'\}$, then $V_D$ is of the form given in the middle in Table 1.
%Applying the substitutions $a=0*$,$a'=1*$,$b=*0$ and $b'=*1$ in words $u\in D$ as described above, we obtain the cube tiling code $V$ (on right in Table 1).

$$
%\begin{equation}
%\label{d1}
%\begin{table}
%\begin{subtable}[c]{0.5\textwidth}
%\centering
\begin{tabular} {c c c c c c c}
 $1$ & $2$ & $3$ & 4 & 5 & 6 & 7\\
\hline
$0$ & $*$ & $0$ & $*$& $0$ & $*$ & $*$\\
$1$ & $*$ & $1$ & $*$& $1$ & $*$ & $*$\\
$*$ & $0$ & $0$ & $*$& $1$ & $*$ & $*$\\
$*$ & $1$ & $0$ & $*$& $1$ & $*$ & $*$\\
$0$ & $*$ & $1$ & $*$& $*$ & $1$ & $*$\\
$0$ & $*$ & $1$ & $*$& $*$ & $0$ & $*$\\
$1$ & $*$ & $*$ & $0$& $0$ & $*$ & $*$\\
$1$ & $*$ & $*$ & $1$& $0$ & $*$ & $*$\\
\end{tabular}
%\end{subtable}
%\begin{subtable}[c]{0.5\textwidth}
%\centering
\qquad
\begin{tabular} {c c c}
 $1$ & $2$ & $3$\\
\hline
$a_1$ & $a_3$ & $a_5$\\
$a_1'$ & $a_3'$ & $a_5'$\\
$a_2$ & $a_3$ & $a_5'$\\
$a_2'$ & $a_3$ & $a_5'$\\
$a_1$ & $a_3'$ & $a_6$\\
$a_1$ & $a_3'$ & $a_6'$\\
$a_1'$ & $a_4$ & $a_5$\\
$a_1'$ & $a_4'$ & $a_5$\\
\end{tabular}
\qquad
\begin{minipage}[c]{0.20\textwidth}
%\vspace{-0mm}
\includegraphics[width=2.5cm]{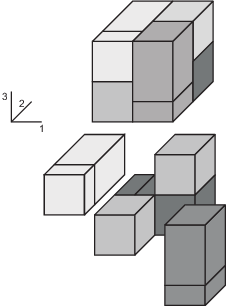}\\
\end{minipage}
%\end{subtable}
%\end{table}
%\end{equation}
$$

\noindent{\footnotesize Table 1: A matrix $M(D)$ of a regular $4$-domino tiling $D$ of the box  $[0,2]^7$ with $*$-column at position $7$ (on the left), its cube tiling code $V_D$ (in the middle) and a realization of $V_D$  (on the right). In this realization we took $X=[0,1]^3, f_1(a_1)=[0,\frac{1}{2}), f_1(a_2)=[0,\frac{1}{4}), f_2(a_3)=[0,\frac{1}{2}), f_2(a_4)=[0,\frac{1}{4})$ and  $f_3(a_5)=[0,\frac{1}{2}), f_3(a_6)=[0,\frac{1}{4})$.}

}
\end{nap}

\medskip
It follows from the above that the meta-structure of an $n$-dimensional regular $(n-d)$-domino tiling $D$ of $[0,2]^n$ is exactly like the structure of the $d$-dimensional cube tiling code $V_D\subset S^d$. More precisely, let $\prop(D)=\bigcup_{v\in D}\prop(v)$, and let $X_i=[0,2]^{n_i}$ for $i\in [d]$, where $n_i=|N_i|$ (the sets $N_i$, $i\in [d]$, are of the form (\ref{roz})). For simplicity of the notation, we take $\prop(D)=[n]$ (that is $D$ is $*$-irreducible) and for every $i\in [d]$ the set $N_i$ consists of consecutive integers (compare Table 1).
Let us consider the following realization of $V_D$: For $v\in V_D$ and $i\in [d]$ if $v_i=a_j$, where $j\in N_i$ ($a_j,a_j'$ are such as in (\ref{rel})), then 
$$
f_i(a_j)=[0,2]\times\cdots \times [0,2]\times [0,1)\times [0,2]\times \cdots \times [0,2], 
$$
where $[0,1)$ stands at the $j$-th position ($f_i(a_j)$ consists of $n_i$ factors) and 
$$
f_i(a_j')=[0,2]\times\cdots \times [0,2]\times [1,2]\times [0,2]\times \cdots \times [0,2],
$$
 where $[1,2]$ stands at the $j$-th position  ($f_i(a_j')$ consists of $n_i$ factors). Thus, $f_i(a_j')=X_i\setminus f_i(a_j)$ for $i\in [d]$ and $j\in N_i$. Let $\ka D=\{f_1(v_1)\times \cdots \times f_d(v_d)\colon v\in V_D\}$. Since $V_D$ is a cube tiling code, its realization $\ka D$ is a minimal partition of the $d$-box $X=X_1\times \cdots \times X_d$ (here we ignore the dimensions of $X_i$ for $i\in [d]$). On the other hand, to pass from $\ka D$ to $D$ we take into account the dimensions of $X_i$ for $i\in [d]$. To do that, let $A\in \{{\rm Cl}(B)\colon B\in \ka D\}$,  where ${\rm Cl}(B)$ is the closure of $B$. Then $A=A(v)=A_1\times\cdots \times A_d$, where $v\in V_D$ and the sets $A_i=A^1_i\times \cdots \times A^{n_i}_i$, $i\in [d]$, are of the form  ${\rm Cl}(f_i(a_j))$ or $f_i(a_j')$. Since $v\in V_D$, there is $u\in D$ such that $v=v(u)$. Then 
$$
\breve{u}=A^1_1\times \cdots \times A^{n_1}_1 \times A^1_2\times \cdots \times A^{n_2}_2 \times \cdots \times A^1_d\times \cdots \times A^{n_d}_d.
$$
%Moreover, for every $A\in \ka D$ and every $((x^1_1,...,x^1_{n_1}),...,(x^d_1,...,x^d_{n_d}))\in A$

\medskip
Clearly, every cube tiling code $V\subset S^d$ induces a regular $(n-d)$-domino tiling of $[0,2]^n$. Indeed, let $V\subset S^d$ be a cube tiling code, and let $S_i\subset S$ be the set of all letters that appear in words of $V$ at the $i$-th position, that is, $S_i=\{v_i\colon v\in V\}$. Since $V$ is a cube tiling code, if $l\in S_i$, then $l'\in S_i$, and therefore $S_i$ is of the form $S_i=\{a^1_i,(a^1_i)',..., a^{r_i}_i,(a^{r_i}_i)'\}$ for $i\in [d]$, $r_i\geq 1$.  Let us recall, that letters of an alphabet $S$ (see Section 1) were arbitrary objects with the property $s''=(s')'=s$ and $s'\neq s$. In particular, we can take words of the form $*...*\delta*..*$ as letters, where $\delta \in \{0,1\}$ and $(*...*0*..*)'=*...*1*..*$ (0,1 stand at the same position). Therefore,  for every $i\in [d]$ we define a new alphabet $\bar{S}_i=\{\delta *\ldots *,\;\; *\delta *\ldots *, ..., *...*\delta \colon \delta \in \{0,1\}\}$, where  $|\bar{S}_i|=2r_i$ for $i\in [d]$ and the length of each word in $\bar{S}_i$ is $r_i$. Thus, $\bar{S}_i$ arises from $S_i$ by replacing each pair $a^j_i,(a^j_i)'$, $j\in [r_i]$, by the pair $*...*0*..*, *...*1*..*$, where $0$ and $1$ stand at the $j$-th position.

\pagebreak
Let now

$$
\bar{v}_i=
\begin{cases}
*...*0*...* & \text{if $v_i=a_i^j$, where $0$ stands at the $j$-th position, $j\in [r_i]$},\\
*...*1*...* & \text{if $v_i=(a_i^j)'$, where $1$ stands at the $j$-th position, $j\in [r_i]$}.\\
\end{cases}
$$

\medskip
It is easy to see that the set $D_V=\{\bar{v}=\bar{v}_1...\bar{v}_d\colon v\in V\}$ is a regular $(n-d)$-domino tiling of $[0,2]^n$ with $K^*=\emptyset$ (that is, $D_V$ is $*$-irreducible). We shall say that $D_V$ is {\it induced} by $V$. Let us note also, that if the code $V$ is written down in a way that $S_i\cap S_j=\emptyset$ for $i\neq j, i,j\in [d]$, then $V$ is the cube tiling code of $D_V$. For example, if $V=\{aa,aa',a'b, a'b'\}$, then $S_1=\{a,a'\}, S_2=\{a,a',b,b'\}$, but $V$ is not a cube tiling of the tiling $D_V$, as $S=\{a,a',b,b'\}$. Such $V$ has to be rewritten in the form  $V=\{a_1a_2,a_1a_2',a_1'a_3,a_1'a_3'\}\subset \{a_1,a_1',a_2,a_2',a_3,a_3'\}^2$ to be (formally) a cube tiling code of an $(n-d)$-domino tiling.

\medskip
At the end of this section we give a characterization of numbers $n,d\in \enn, n>d$, for which regular $*$-irreducible $(n-d)$-domino tilings of the box $[0,2]^n$ exist. To do that, we need the following property of cube tiling codes which is a consequence of  \cite[Theorem 6]{Pr}.

%we back to the problem of the existence of  a regular $*$-irreducible $(n-d)$-domino tiling of the box $[0,2]^n$ for a given $n,d\in \enn, n>d$. We shall need the following property of cube tiling codes which is a consequence of  \cite[Theorem 6]{Pr}.

\begin{lemat}
\label{lite}
Let $d\in \enn$. For every $n\in \{d,...,2^d-1\}$ there is a cube tiling code $V\subset S^d$ such that $\sum_{i=1}^dl_i=n$, where $l_i=|\{v_i\colon v\in V\}|/2$ for $i\in [d]$.
\end{lemat}
\proof
Let $W\subset S^d$ be a cube tiling code, and let $m_i=|\{w_i\colon w\in W\}|/2$ for $i\in [d]$. For $l\in S$ and $i\in [d]$ let $W^{i,l}=\{w\in W\colon w_i=l\}$. If $\sum_{i=1}^dm_i>d$, then there are $a,b\in S$ such $\{a,a'\}\cap \{b,b'\}=\emptyset$ and the sets $W^{i,a}\cup W^{i,a'}$ and $W^{i,b}\cup W^{i,b'}$ are non-empty (and clearly disjoint). Note that, since $W$ is a code, for every $w\in W^{i,a}\cup W^{i,a'}$ and $u\in W^{i,b}\cup W^{i,b'}$ there is $j\in [d]\setminus \{i\}$ such that $w_j=u_j'$. Therefore, the set $W^1\subset S^d$ defined by the formula
$$
W^1=W\setminus (W^{i,b}\cup W^{i,b'})\cup \{w_1...w_{i-1}aw_{i+1}...w_d\colon w \in W^{i,b}\}\cup \{w_1...w_{i-1}a'w_{i+1}...w_d\colon w \in W^{i,b'}\}
$$
is a cube tiling code, and moreover $\sum_{i=1}^dm^1_i=\sum_{i=1}^dm_i-1$, where $m^1_i=|\{w_i\colon w\in W^1\}|/2$. Clearly, in this way one may show that if $\sum_{i=1}^dm_i=d+k$ for some $k\in \enn$, then for every $r\in [k]$ there is a cube tiling code $W^r$ such that $\sum_{i=1}^dm^r_i=\sum_{i=1}^dm_i-r$, where $m^r_i=|\{w_i\colon w\in W^r\}|/2$.

Now our result follows from \cite[Theorem 6 for $n=2$]{Pr} which implies that  $\sum_{i=1}^dl_i\leq 2^d-1$ for every cube tiling $V\subset S^d$, where $l_i=|\{v_i\colon v\in V\}|/2$ and there is a cube tiling code $U\subset S^d$ such that  $\sum_{i=1}^dl_i=2^d-1$, where $l_i=|\{u_i\colon u\in U\}|/2$. 
\hfill{$\square$}

\begin{uw}{\rm
In \cite[Theorem 6]{Pr} a cube tiling $\ka T=[0,1)^d\oplus T$ of the $d$-dimensional torus $\mathbb{T}^d_m$ is considered. Since $\ka T$ is a minimal partition, to get the inequality  $\sum_{i=1}^dl_i\leq 2^d-1$ for a cube tiling code $V\subset S^d$, we just treat the tiling $\ka T$ as a realization of $V$  (see also \cite[Theorem 7]{Pr}). 

In the last step of the proof of Lemma \ref{lite} instead of \cite[Theorem 6]{Pr} we may use \cite[Corollary]{Als}, however, this requires a transition from a cube tiling code $V\subset S^d$ to the regular $(n-d)$-domino tiling $D_V$ of $[0,2]^n$. 
}
\end{uw}

\begin{tw}
\label{lete}
Let $n,d\in \enn$ and $n>d$. There is a regular $*$-irreducible $(n-d)$-domino tiling of the box $[0,2]^n$ if and only if $n\in \{d+1,...2^d-1\}$.
\end{tw}
\proof
Let $D$ be a regular $*$-irreducible $(n-d)$-domino tiling of the box $[0,2]^n$, and let $V_D$ be the cube tiling code of $D$. For $i \in [d]$ let $l_i=|\{v_i\colon v\in V_D\}|/2$. Note that, for every $i\in [d]$ each word in the set $K_i$  contains $l_i-1$ stars (compare (\ref{ki}) and Table 1). Since $D$ is $*$-irreducible, we have $\sum_{i=1}^d(l_i-1)=n-d$, and thus, $\sum_{i=1}^dl_i=n$. From \cite[Theorem 6]{Pr} it follows that $n\leq 2^d-1$, and since $n>d$, we have $n\in \{d+1,...,2^d-1\}$. 

Let now $n\in \{d+1,...2^d-1\}$. By Lemma \ref{lete}, there is a cube tiling $V\subset S^d$ such that $\sum_{i=1}^dl_i=n$, where $l_i=|\{v_i\colon v\in V\}|/2$ for $i\in [d]$. Therefore, the tiling $D_V$ induced by $V$ is a regular $*$-irreducible $(n-d)$-domino tiling of the box $[0,2]^n$.
\hfill{$\square$}

\subsection{Flip-connected classes for regular $(n-d)$-domino tilings}

In the proof of our version of Theorem \ref{glue} for regular $(n-d)$-domino tilings we shall need the following property: %(It is rather obvious however we give a proof of it.)

\begin{lemat}
\label{lit}
Let $V,W\subset S^d$ be two cube tiling codes such that one can pass by shifts from $V$ to $W$, and let $S_i\subset S$ be the set of all letters that appear in words of the set $V\cup W$ at the $i$-th position. Then it is possible to pass by shifts $\{v^1,w^1\},...,\{v^k,w^k\}$ from $V$ to $W$ such that $v^j_i,w^j_i\in S_i$ for every $j\in [k]$.   
\end{lemat}
\proof
Assume that there is a sequence of shifts $\{r^1,p^1\},...,\{r^k,p^k\}$ which transforms $V$ into $W$ such that $ r^t_i=l\not\in S_i$ for some $t\in [k]$. As it was noted in Section 1,  there is a sequence $\{u^1,q^1\},...,\{u^k,q^k\}$ of twin pairs and a sequence $V^1,...,V^{k+1}$ of cube tiling codes such that $r^j,p^j$ is a shift of $u^j,q^j$  for $j\in [k]$, and moreover 
$$
V^1=V, \;\;\;\;  V^{j+1}=(V^{j}\setminus \{u^j,q^j\})\cup \{r^j,p^j\}\;\; {\rm and}\;\; V^{k+1}=W,
$$
where $\{u^j,q^j\}\subset V^{j}$ for $j\in [k]$. Let 
$$
(V^j)^{i,l}=\{v\in V^j\colon v_i=l\}\;\; {\rm and}\;\; (V^j)^{i,l'}=\{v\in V^j\colon v_i=l'\}.
$$
For $s,s'\in S\setminus \{l,l'\}$ let 
$$
U^j=(V^j\setminus ((V^j)^{i,l}\cup (V^j)^{i,l'}))\cup (Q^j)^{i,s}\cup (Q^j)^{i,s'}, 
$$
where 
$$
(Q^j)^{i,s}=\{v_1...v_{i-1}sv_{i+1}...v_d\colon v\in (V^j)^{i,l}\}\;\; {\rm and}\;\; (Q^j)^{i,s'}=\{v_1...v_{i-1}s'v_{i+1}...v_d\colon v\in (V^j)^{i,l'}\}.
$$
It is easy to check that $U^j$ is a cube tiling code for $j\in [k+1]$ (compare \cite[Subsection 2.1]{Kis}). Since $l,l'\not\in S_i$, we have  $U^1=V$ and $U^{k+1}=W$. 

Let  now $(\bar{u}^j,\bar{q}^j)_{j\in [k]}$ and $(\bar{r}^j,\bar{p}^j)_{j\in [k]}$ be two sequences of twin pairs, which arise from $(u^j,q^j)_{j\in [k]}$ and $(r^j,q^j)_{j\in [k]}$ by the substitutions $l\rightarrow s, l'\rightarrow s'$ at the $i$-th position in every word of  $(u^j,q^j)_{j\in [k]}$ and $(r^j,q^j)_{j\in [k]}$.  By the definition of $U^j$, we have $\{\bar{u}^j,\bar{q}^j\}\subset U^{j}$ and $U^{j+1}=(U^{j}\setminus \{\bar{u}^j,\bar{q}^j\})\cup \{\bar{r}^j,\bar{p}^j\}$ for $j\in [k]$.

%By the definition of $U^j$, two sequences $(\bar{u}^j,\bar{q}^j)_{j\in [k]}$ and $(\bar{r}^j,\bar{p}^j)_{j\in [k]}$ of twin pairs, which arise from $(u^j,q^j)_{j\in [k]}$ and $(r^j,q^j)_{j\in [k]}$ by the substitutions $l\rightarrow s, l'\rightarrow s'$ at the $i$-th position in every word of  $(u^j,q^j)_{j\in [k]}$ and $(r^j,q^j)_{j\in [k]}$ are such that $U^{j+1}=(U^{j}\setminus \{\bar{u}^j,\bar{q}^j\})\cup \{\bar{r}^j,\bar{p}^j\}$, where $\{\bar{u}^j,\bar{q}^j\}\subset U^{j}$ for $j\in [k]$. 

Since $r^j,p^j$ was a shift of $u^j,q^j$  for $j\in [k]$, it follows that $\bar{r}^j,\bar{p}^j$ is a shift of $\bar{u}^j,\bar{q}^j$  for $j\in [k]$. Thus, we may pass from $V$ to $W$ by the sequence of shifts $\{\bar{r}^1,\bar{p}^1\},...,\{\bar{r}^k,\bar{p}^k\}$. As there was no additional assumption on the letter $s\in S$, we may assume that $s,s'\in S_i$.
\hfill{$\square$}

%the code $(U^j)^{i,s}\cup (U^j)^{i,s'}$ arises from $(V^j)^{i,l}\cup (V^j)^{i,l'}$ by changing at the $j$-th positions in every word of $(V^j)^{i,l}\cup (V^j)^{i,l'}$ the letters $l,l'$ into the letters $s,s'$.

%Define a new sequence $(U^m)_{m=0}^k$ of cube tiling codes

%such that $v^j,w^j$ is a shift of $u^j,q^j$  for $j\in [k]$, $V^0=V$, $V^j=(V^{j-1}\setminus \{u^j,q^j\})\cup \{v^j,w^j\}$, where $\{u^j,q^j\}\subset V^{j-1}$ for $j\in [k]$ and $W=V^k$.
%Let us note that if there is $r\neq j$ such $u^j=u^r$ and $\{u^{m},q^{m}\}$ is a shift of  $\{u^r,q^r\}$ in the $p$-th direction, where $p\neq i$, then $u_i^{r}=q_i^{r}=u^{m}_i=q^{m}_i=l$. Clearly, since $l\not \in S_i$, at some stage of the process of passing from $V$ to $W$ the words with $l$ at the $i$-th position have to vanish, that is, their shifts contain at the $i$-th position the letters from $S_i$. Thus, we may repeat the whole passing from $V$ to $W$ replacing at the $i$-th position in  $\{u^1,q^1\},...,\{u^k,q^k\}$ the letter $l$ with any other $l^1\in S$. In particular, we may take $l^1\in S_i$.
%\hfill{$\square$} 

\medskip
Let $D$ be an $(n-d)$-domino tiling of $[0,2]^n$. A subset $C\subset D$ is called a {\it simple component} of $D$ if $\prop(v)=\prop(w)$ for every $v,w\in C$ and $\prop(v)\neq \prop(u)$ for every $v\in C$ and $u\in D\setminus C$ (the tiling in Table 1 has four simple components). If  $\prop(v)=\prop(w)$ for every $v,w\in D$, then we say that $D$ is a {\it simple} $(n-d)$-domino tiling of $[0,2]^n$. It is obvious that every simple $(n-d)$-domino tiling of $[0,2]^n$ is regular. Moreover, if $W\subset S^d$ is the cube tiling code of a regular $(n-d)$-domino tiling $D$ of $[0,2]^n$ and $W$ is simple, then $D$ is simple.

Theorem \ref{glue} for regular $(n-d)$-domino tilings of $[0,2]^n$ reads as follows:

\begin{tw}
\label{fl1}
Let $d\leq 6$. The family of all regular $(n-d)$-domino tilings of the box $[0,2]^n$ is flip-connected.
\end{tw}
\proof
Let $D,E\subset (*A)^n$ be two regular $(n-d)$-domino tilings of $[0,2]^n$. 

Let $V_D\subset S^d$ be the  cube tiling code of $D$, $C\subset D$ be a simple component of $D$, and let $V_C\subset V_D$ be the set of all words which encode $C$. By Theorem \ref{glue} we can pass by shifts from $V_D$ to $W$, where $W$ is the simple cube tiling code such that $V_C\subset W$. By Lemma \ref{lit} we may do that using the letters from the set $S_i\subset S$, $i\in [d]$,  which consists of all letters that appear in words in $V_D$ at the $i$-th position (we have $V_C\cup V_D=V_D$ as $V_C\subset V_D$). Let $\{v^1,w^1\},...,\{v^k,w^k\}$ be a sequence of shifts that transform $V_D$ into $W$. 

Let $v,w$ be a twin pair in $V_D$ such that $v_i=w_i'$, and let the twin pair $v^1,w^1$ be the shift of $v,w$, that is,  $v^1_i=(w^1_i)'$, $\{v_i,w_i\}\neq \{v^1_i,w^1_i\}$ and $v_j=w_j=v^1_j=w^1_j$ for $j\in [d]\setminus \{i\}$. We show that the shift $\{v,w\}\rightarrow \{v^1,w^1\}$ induces a flip in $D$ which transform $D$ into a regular $(n-d)$-domino tiling of $[0,2]^n$. To do that, let us recall that $S=\{a_1,a_1',...,a_n,a_n'\}$. Thus,  $v_i,w_i=a_j,a_j'$ and $v^1_i,w^1_i=a_s,a_s'$, where $j\neq s$ and $j\in N_i$, where $N_i$ is as in (\ref{roz}) (clearly,  $j\in N_i$ by the definition of $N_i$). Note now, that, by our assumption stemming from Lemma \ref{lit}, we have $s \in N_i$. Let $e,f\in D$ be the twin pair such that the words $v,w$ encode it. We flip $e,f$ into $e^1,f^1$, where $e^1_j=f^1_j=*$, $\{e^1_s,f^1_s\}=\{0,1\}$ and $e_m=f_m=e^1_m=f^1_m$ for $m\in [n]\setminus\{j,s\}$.  Let $D^1=(D\setminus \{e,f\})\cup \{e^1,f^1\}$. Clearly, if $V_D$ contains at least two words with the letter $a_j$, then partitions of the form  (\ref{roz}) for $D$ and $D^1$ are the same. Thus, $D^1$ is regular. Moreover, $V_{D^1}=(V\setminus \{v,w\})\cup \{v^1,w^1\}$ is the cube tiling code of $D^1$. 

Suppose now that $v$ is the only word in $V_D$ with $v_i=a_j$, that is, the $j$-th column in $M(D)$ contains only one pair $0,1$ ($0$ for $v_i$ and $1$ for $w_i=v_i'$). Let $[n]=N_1\dot{\cup} \cdots \dot{\cup} N_d \dot{\cup} N^*$ be as in (\ref{roz}) for the tiling $D$. After gluing the boxes $e,f$ at the $j$-th position that pair $0,1$ vanishes and the $j$-th column in $M(D^1)$ consists of stars only. This means that that the partition $[n]=M_1\dot{\cup} \cdots \dot{\cup} M_d \dot{\cup} M^*$, where $M^*=N^*\cup \{j\}$, $M_m=N_m$ for $m\neq i$ and $M_i=N_i\setminus \{j\}$, has the properties ($\alpha$) and ($\beta$) for the tiling $D^1$. To see this, observe that $M_i\neq\emptyset$, as $|N_i|\geq 2$. Moreover, the set $\bar{K}_i$ of the rows of the sub-matrix $[d_{ij}]_{i\in [2^d], j\in M_i}$ of the matrix $M(D^1)$ is of the form (\ref{ki}): It has two elements less than the set $K_i$ for the tiling $D$, and the lengths of words in $\bar{K}_i$ are one less than the lengths of words in $K_i$. The set $M^*$ is non-empty and contains one word of length $|N^*|+1$ consisting only of stars.

%Let $K_1$ ,...,$K_d$ be sets (\ref{ki}) for $D$. Let us note, if $|N_i|=2$, then the set $\bar{K}_i$ is equal to $\{0,1\}$,  and if $|N_i|>2$, then the set $\bar{K}_i$ has two element less than $K_i$, where $K_i$ is of the form (\ref{ki}) for $D$. Hence, $D^1$ is regular. 

Therefore, the sequence of shifts $\{v^1,w^1\},...,\{v^k,w^k\}$ induces the sequence $D^1,...,D^{k+1}$ of regular $(n-d)$-domino tilings of the box $[0,2]^n$ such that $D^{i+1}$ is obtained from $D^i$ by a flip for $i\in [k]$ and $D^{k+1}$ is the simple $(n-d)$-domino tiling of $[0,2]^n$ whose cube tiling code is $W$. Now a procedure of passing from $D$ to $E$ by flips is clear: We pass by flips from $D$ to $C_D$ and $E$ to $C_E$, where $C_D$ and $C_E$ are simple  $(n-d)$-domino tilings of $[0,2]^n$, and next we pass by flips from $C_D$ to $C_E$ (the fact that such passing is possible is obvious). 
\hfill{$\square$}

\medskip
Let $D,E\subset (*A)^n$ be $(n-d)$-domino tilings of $[0,2]^n$, and let $D\approx E$ if and only if one can pass by flips from $D$ to $E$. Clearly, $\approx$ is an equivalent relation on the family of all $(n-d)$-domino tilings of $[0,2]^n$. Let $[D]_\approx$ be the equivalent class of $D$. Note that, if $D$ is regular, then some elements of $[D]_\approx$ need not to be regular. For example, if in the tiling given in Table 1 we flip the twin pair $v^3=*00*1**$ and $v^4=*10*1**$ into $u=**001**$ and $q=**011**$, then we obtain a tiling which is not regular. (That is why in the proof of Theorem \ref{fl1} we needed Lemma \ref{lit}.) We shall call such non-regular elements of  $[D]_\approx$, where $D$ is regular, {\it almost} regular. 

From Theorem \ref{fl1} we obtain

\begin{wn}
\label{alm}
Let $d\leq 6$. If $D$ and $E$ are regular $(n-d)$-domino tilings of the box $[0,2]^n$, then $[D]_\approx=[E]_\approx$. In other words, the family of all almost regular $(n-d)$-domino tilings of the box $[0,2]^n$ is flip-connected.
\end{wn}
\vspace{-3mm}
\hfill{$\square$}

\medskip
At the end of the paper we give an example of a class of  subsets of the box $[0,2]^n$ with the property that every family of regular $(n-d)$-domino tilings from that class is flip-connected. To do that, we need the following  result which stems from \cite[Theorem 5.1]{Kgl}:

\begin{tw}
\label{fl4}
Let $d\leq 5$, $D$ be a  regular $(n-d)$-domino tiling of $[0,2]^n$ with a simple component $C\subset D$, and let $Q$ be the simple $(n-d)$-domino tiling of $[0,2]^n$ such that $C\subset Q$. It is possible to pass by flips from $D$ to $Q$ keeping $C$ fixed.
\end{tw}
\proof
Let $V_D,V_C$ and $W$ be as in the proof of Theorem \ref{fl1}.
%Let $V_D$ be the cube tiling code of $D$, $V_C\subset V_D$ be the simple component of $V_D$ that encoded $C$, and let $W$ be the simple cube tiling code such that $V_C\subset W$.  
By \cite[Theorem 5.1]{Kgl} there is a sequence $\{v^1,w^1\},...,\{v^k,w^k\}$ of shifts that transform $V_D$ into $W$ keeping $V_C$ fixed.
%that is, the sets $V_C$ and $\{v^1,w^1,\\...,v^k,w^k\}$ are disjoint. 
Precisely in the same manner as in the proof of Theorem \ref{fl1} we show that the sequence of shifts $\{v^1,w^1\},...,\{v^k,w^k\}$ induces the sequence $D^1,...,D^{k+1}$ of regular $(n-d)$-domino tilings of the box $[0,2]^n$ such that $D^{i+1}$ is obtained from $D^i$ by a flip for $i\in [k]$ and $D^{k+1}$ is the simple $(n-d)$-domino tiling of $[0,2]^n$ whose cube tiling code is $W$. Clearly, $D^{k+1}=Q$ and since $V_C$ was unchanged during the passing from $V_D$ to $W$, the set $C\subset Q$ was unchanged during the passing from $D$ to $Q$.
\hfill{$\square$}

\medskip
Let $C$ be a simple component of  a regular $(n-d)$-domino tiling $D$ of $[0,2]^n$, and let $F_C=[0,2]^n\setminus \bigcup_{v\in C} \breve{v}$. Let us note, that $F_C$ does not have to be a box. In what follows we describe a family of $(n-d)$-domino tilings  of $F_C$ which is flip-connected.  
%From \cite[Corollary 5.2]{Kgl} we obtain

\begin{wn}
\label{fl5}
Let $d\leq 5$, $C$ be a simple component of  a regular $(n-d)$-domino tiling $D$ of $[0,2]^n$, and let and let $F_C=[0,2]^n\setminus \bigcup_{v\in C} \breve{v}$. Then the family $\ka F$ of all $(n-d)$-domino tilings of the set $F_C$ such that for every $F\in \ka F$ the tiling $F\cup C$ is regular is flip-connected.
\end{wn}
\proof
Let $F,G\in \ka F$, and let $Q$ be the simple $(n-d)$-tiling of $[0,2]^n$ such that $C\subset Q$. Since $F\cup C$ and $G\cup C$ are regular, by Theorem \ref{fl4} there are two sequences $F^1,...,F^k\in \ka F$ and $G^1,...,G^m\in \ka F$ such that $F=F^1, Q=F^k$ and  $G=G^1, Q=G^m$ and $F^{i+1}$ is obtained from $F^i$ by a flip for $i\in [k-1]$, and similarly $G^{i+1}$ is obtained from $G^i$ by a flip for $i\in [m-1]$. Thus, $G$ can be obtained from $F$ by a sequence of flips.%$F^1,...,F^k,G^{m-1},...,G^1$.
\hfill{$\square$}

\begin{uw}
{\rm In the definition of $F_C$ we used the tiling $D$ but we can omit it, leaving only simple component $C$. To do that, let $Q$ be a simple $(n-d)$-domino tiling of $[0,2]^n$, where sets (\ref{roz}) are of the form $N_1=\{i_1\},...,N_d=\{i_d\}$. For $i\in \{i_1,...,i_d\}$ let $P_i$ be a family of twin pairs $A,B\in Q$ such that $\{A_i,B_i\}=\{[0,1],[1,2]\}$. Using  \cite[Theorem 1]{kpcoi} it can be shown that if $C=Q\setminus \bigcup_{i\in M} P_i$, where $M\subset \{i_1,...,i_d\}$, then
$F_C=\bigcup_{v\in P}\breve{v}$, where $P=\bigcup_{i\in M} P_i$.

% \cite[Theorem 1]{kpcoi} it follows that if $D$ is a regular $(n-d)$-domino tiling with a simple component $C\subset Q$, then $C=Q\setminus \bigcup_{i\in M} P_i$, where $M\subset \{i_1,...,i_d\}$ and $P_i$, $i\in M$, is  a family of twin pairs $A,B\in Q$ such that $\{A_i,B_i\}=\{[0,1],[1,2]\}$.

%can be shown that if $C=Q\setminus \bigcup_{i\in M} P_i$, where $M\subset \{i_1,...,i_d\}$, then
%$F_C=[0,2]^n\setminus \bigcup_{v\in P}\breve{v}$, where $P=\bigcup_{i\in M} P_i$.
}
\end{uw}
%The definition of the set $F_C$ needs  the tiling $D$. However, there is a nice characterization of sets like $F_C$ (\cite{kpcoi}) which allows us reformulate the above corollary such that $D$ is not necessary.

%Two codes $V,W\subset (*A)^n$ are {\it isomorphic} if there are permutations of columns and rows and there are mappings $l\mapsto l'$ of letters in columns of the matrix $M(V)$, where $0'=1, 1'=0$ and $*'=*$,  which transform $M(V)$ into $M(W)$.

%A word $q=q_1...q_d$ is an {\it edge} in $\{0,1\}^d$ if $q_i=*$ for exactly one $i\in [d]$ and $q_j\in \{0,1\}$ every $j\in [d]\setminus \{i\}$, and two edges $q,p$ are {\it independent} if they are dichotomous (that is, $\{q_i,p_i\}=\{0,1\}$ for some $i\in [d]$). 

%From \cite[Theorem 1]{kpcoi} it follows that  Corollary \ref{fl5} can be given in the form:

%\begin{wn}
%\label{fl6}
%Let $d\leq 5$, and let $W=\{q^i*^{n-d-1}\colon i\in [k]\}$, where $q^1,...,q^k$ are pairwise independent edges in $\{0,1\}^d$. For every code $Q$ which is isomorphic to $W$, the family of all almost regular $(n-d)$-domino tilings of the set $\bigcup_{v\in Q}\breve{v}\subseteq [0,2]^n$ is flip-connected.
%\hfill{$\square$}
%\end{wn}

\section{Final comments}

%\medskip
%\noindent{{\sc Remarks}
Since Keller's cube tiling conjecture in dimension $d=7$ is true (\cite{BHM,Kis}), every almost regular $(n-7)$-domino tiling of $[0,2]^n$ has a twin pair. However, we do not know whether Theorem \ref{fl1} for $d=7$ is true. The method applied in \cite{Kgl} (which is based on computations) seems to be non-adequate for the case $d=7$ because of a computational complexity. For all dimensions $d\geq 8$ Theorem \ref{fl1} is not true (\cite{LS2,M}): If $M$ is the twin pair free cube tiling code in dimension eight (\cite{M}), then the regular $8$-domino tiling $D_M$ of the box $[0,2]^{16}$ whose cube tiling code is $M$ does not contain a twin pair, and thus, $|[D_M]_\approx|=1$. (More about cube tiling codes without twin pairs can be found in \cite{LS2}.)  In general, one can find a twin pair free $(n-d)$-domino tiling of $[0,2]^n$ for $d\geq 3$. For example, $D=\{*01\varepsilon,1*0\varepsilon,01*\varepsilon\colon \varepsilon \in \{0,1\}\}\cup \{000*,111*\}$ is a domino tiling of $[0,2]^4$ without twin pairs (that is, $|[D]_\approx|=1$).

The family of all cube tiling codes $V\subset S^d$ is known for $d\leq 4$: In \cite{Kii} we described the family $\ka N_d$ of all non-isomorphic cube tiling codes for $d=4$. We obtained  $|\ka N_4|=27358$ and the number of all cube tiling codes for $d=4$ is equal to 17, 794, 836, 080, 455, 680. For $d=5$ we know that $|\ka N^4_5|=899, 710, 227$, where $\ka N^4_5$ is the family of all non-isomorphic cube tiling codes $V\subset S^5$ with $|S|=4$ (\cite{Ost}). It is known also (\cite{Kgl}) that the cardinality of the family of all cube tiling codes $V\subset S^6$ far exceeds the number $10^{84}$. However, we do not know what a portion of all $(n-d)$-domino tilings of $[0,2]^n$ are almost regular tilings.%}
 %It seems that the class of almost regular $(n-d)$-domino tilings of the box $[0,2]^n$ for $d\leq 6$ is not very small. 

\medskip

%Let $D\subset (*A)^n$ be an $(n-d)$-domino tiling of $[0,2]^n$, and let $\prop(D)=\bigcup_{v\in D}\prop(v)$. Let us call the number $d_e(D)=|\prop(D)|$ the {\it essential dimension} of $D$. Clearly, $d_e(D)\leq n$. In [] (see also []) it was shown that if $D\subset (*A)^n$ is an $(n-d)$-domino tiling of $[0,2]^n$, then $d_e(D)\leq 2^d-1$ and for every $m\in \{d,..., 2^{d-1}\}$ there is a  $(n-d)$-domino tiling $D$ such that $d_e(D)=m$.

\end{document}